\documentclass[12pt]{article}
\usepackage{amsfonts, amsmath}
\usepackage[pdftex]{graphics}
\usepackage[hang]{subfigure}
\usepackage{rotating}

\def\complaint#1{}
\def\withcomplaints{
\newcounter{mycomplaints}
\def\complaint##1{\refstepcounter{mycomplaints}%
\ifhmode%
\unskip%
{\dimen1=\baselineskip \divide\dimen1 by 2 %
\raise\dimen1\llap{\tiny -\themycomplaints-}}\fi%
\marginpar{\tiny [\themycomplaints]: ##1}}%
}

\withcomplaints 
\newtheorem{theorem}{Theorem}
\newtheorem{conj}{Conjecture}

\newtheorem{prop}[theorem]{{Proposition}}

\def\T{\mathbb T}
\def\C{\mathbb C}
\def\g{{g}}
\def\E{\mathbb E}

\def\Z{\mathbb Z}

\def\p{{\bf p}}

\begin{document}

\title {Periodic Planar Disk Packings}

\author{Robert Connelly \thanks{Research supported~ in~ part~ by~ NSF~ Grant~ No. DMS--0209595 (USA). 
\newline  e-mail:
connelly@math.cornell.edu} \\
Department of Mathematics, Cornell University\\
Ithaca, NY 14853, USA\\ \smallskip\\
William Dickinson\thanks{e-mail: dickinsw@gvsu.edu}\\ Department of Mathematics, Grand Valley State University\\ A-2-178 Mackinac Hall, 1 Campus Drive, \\
Allendale, MI 49401-6495, USA}
\maketitle

\begin{abstract}  Several conditions are given when a packing of equal disks in a
torus is locally maximally dense, where the torus is defined as the
quotient of the plane by a two-dimensional lattice. Conjectures are
presented that claim that the density of any collectively jammed packing, whose graph does not consist of all triangles and the torus lattice is the standard triangular lattice, is at most $\frac{n}{n+1}\frac{\pi}{\sqrt{12}}$, where $n$ is the number of packing disks in the torus.  Several classes of collectively jammed packings are presented where the conjecture holds.


\end{abstract}
\section{Introduction and definitions} \label{section:introduction}

Packings of disks are of great importance in number theory, granular materials,  algebraic number theory, and who knows what else.  Their rigidity properties are also important and can provide useful tools for their analysis, in particular for finding maximally dense packings, or locally maximally dense packings, in certain situations.  Here we primarily investigate the case of periodic packings of equal disks, where the period lattice is the ubiquitous triangle lattice given by the edge-to-edge tiling of equilateral triangles.  One reason for choosing this lattice is that it is the basis for the most dense packing of equal disks in the plane.  See, for example, the classic books \cite{L-Fejes-Toth, L-Fejes-Toth2}, for a good discussion of these problems as well as \cite{Grunbaum-tilings} for many examples of square and triangle tilings.  Another reason for this choice is that the techniques of rigidity theory can be applied here, while in many other cases some of the tools described here are not available.

\subsection{Lattices} \label{subsection:ilattices}

A lattice $\Lambda$ in Euclidean space $\E^d$ is the set of all integral linear combinations of the vectors $\g_1, \dots, \g_D$ in $\E^d$,
\begin{equation}\Lambda(\g_1, \dots, \g_D) = \{n_1\g_1 + \dots + n_D\g_n \mid n_1, \dots, n_D \in \Z \}. \label{eqn:lattice-def}
\end{equation}
Note that $D$ can be greater or less than $d$, but the usual definition of a lattice is when $D=d$ and $\g_1, \dots, \g_D$ is a basis of $\E^d$.  This sort of lattice can be regarded as a projection, or injection, of the usual lattice into $\E^d$.  

\begin{figure}[here]
    \begin{center}
        \includegraphics[width=.8\textwidth]{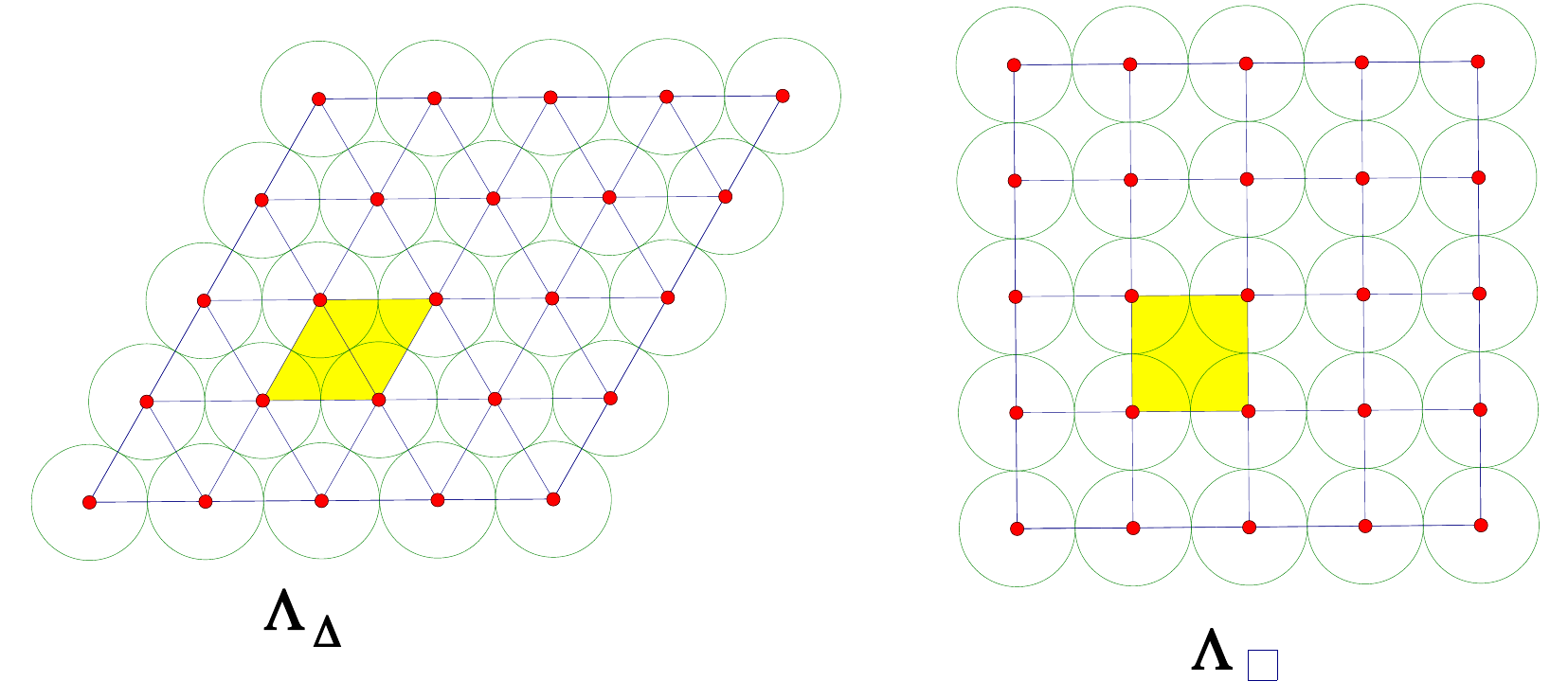}%
        \end{center}
    \caption{This shows a portion of a triangular lattice $\Lambda_{\Delta}$ and the square lattice $ \Lambda_{\Box}$,  along with the corresponding circle lattice packings.  A fundamental region is shaded.  The circle lattice packing for $\T^2(\Lambda_{\Delta})$ has the densest packing density of any packing in the plane.}
    \label{fig:lattices}
    \end{figure}

We also find it is useful, for $d=2$, to use complex notation, so $\C = \E^2$, and define $\g(\theta)=e^{i\theta}=\cos\theta + i \sin \theta$.  For notational fluidity, we define 
\begin{eqnarray*}
\g_{\Delta} = \g(2\pi/6)= \frac{1}{2} + \frac{\sqrt{3}}{2}i \,\, \text{and} \,\,\Lambda_{\Delta} = \Lambda(1, \g_{\Delta}),\\
\g_{\Box} = \g(2\pi/4) = i  \,\, \text{and} \,\,\Lambda_{\Box} = \Lambda(1, \g_{\Box}),
\end{eqnarray*}
which is the usual triangular lattice and the square lattice as in Figure \ref{fig:lattices}.  We say that any lattice of the form $\Lambda(z, g_{\Delta} z)$ is a \emph{triangular lattice}.  In other words there are generators of equal length such that the angle between them is $2\pi/6$.  Similarly a \emph{square lattice} is of the form $\Lambda(z, g_{\Box} z)$.

For any two-dimensional lattice $\Lambda$ in the plane, not contained in a line, we define the corresponding torus as
\begin{equation*}
\T^2(\Lambda) = \C/\Lambda,
\end{equation*}
where two points are identified if they differ by a vector in the lattice $\Lambda$.   We regard $\T^2(\Lambda)$ as a compact metric space, where locally it is isometric to the flat Euclidean plane.  Alternatively, we can regard $\T^2(\Lambda)$ as obtained by identifying opposite edges of a fundamental parallelogram as in Figure \ref{fig:lattices}.  

Since  $\T^2(\Lambda)$ is compact, it has a finite area, namely the area of any fundamental parallelogram.  For example, 
\begin{equation}\label{eqn:area}
\text{Area}(\T^2(\Lambda(\g_1, \g_2))=|\g_1| |\g_2|\sin \theta = \text{Im}(\g_1\bar{\g}_2),
\end{equation}
where $\theta$ is the angle between $\g_1$ and $\g_2$, $\bar{z}$ is the complex conjugate of $z$, and Im($z$) is the imaginary part of $z$.

We say that a \emph{triangular torus} is any torus of the form $\T^2(\Lambda(z, \g_{\Delta}z))$.  In other words, a triangular torus is where the lattice used to define it is a triangular lattice.  Similarly a \emph{square torus} is defined by a square lattice.
    
\subsection{Packings} \label{subsection:packings}

A collection of $n$ equal (radii) disks with disjoint interiors  in $\T^2(\Lambda)$ is called a \emph{packing}.   When $n=1$, the packing is called a \emph{lattice packing}. The density of a packing is $\delta$, the ratio of the sum of the areas of the disks divided by the area of the torus.  So the density of $n$ disks of radius $r$ in $\T^2(\Lambda(z, \g(\theta)z))$ is

\begin{equation}
\delta =n \pi r^2/(|z|^2\sin \theta).\label{eqn:gen-density}
\end{equation}

For the triangular lattice packing when $r=1/2$ and $n=1$, 

\begin{equation}
\delta_{\Delta} = \pi (\frac{1}{2})^2/\sin (2\pi/6) = \pi/\sqrt{12},\label{eqn:max-density}
\end{equation}
which is the maximum density for any packing of equal disks in the plane.

When the packing is lifted to a packing of an infinite number of packing disks in the plane, this definition of packing density agrees with the usual one, where larger and larger regions are intersected with the packing to calculate the packing density.

\begin{figure}[here]
    \begin{center}
        \includegraphics[width=.4\textwidth]{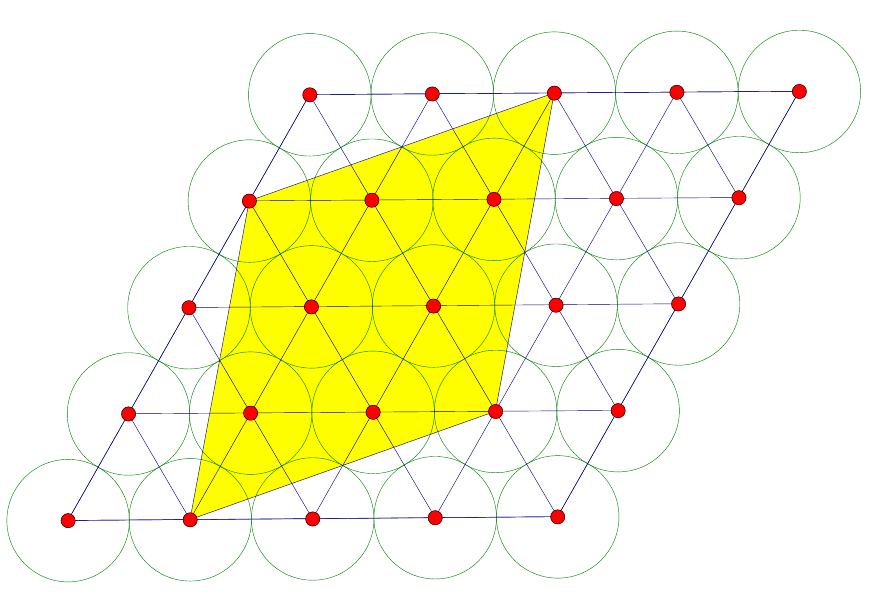}%
        \end{center}
    \caption{This shows a triangular packing of $7$ equal disks in a triangular torus of maximal density $\pi/\sqrt{12}$.  The fundamental region of the sublattice is indicated.}
    \label{fig:triangular-packing}
    \end{figure}

One of our main problems is to understand the most dense packings of $n$ equal disks in a triangular torus.  For certain values of $n$, this determination is easy, namely when there is a triangular lattice packing inside the triangular torus.  This happens when $\Lambda_{\Delta}$ is a lattice containing $\Lambda(z, g_{\Delta}z)$.  So $z = n_1 + n_2g_{\Delta}$, and 
\begin{eqnarray}
|z|^2&=&z\bar{z}=(n_1 + n_2g_{\Delta})(n_1 + n_2\bar{g}_{\Delta})=n_1^2 + n_1n_2(g_{\Delta}+\bar{g}_{\Delta}) +n_2^2 \nonumber \\
&=& n_1^2 + n_1n_2 + n_2^2.\label{eqn:triangle-number}
\end{eqnarray}
The ratio of the areas of the tori is given by (\ref{eqn:triangle-number}), since the angle between the generators $2\pi/6$ is the same in both cases.  Furthermore, (\ref{eqn:triangle-number}) is the number of packing disks by (\ref{eqn:max-density}) since the density is the same for both packings.  Putting this together, we get the following:

\begin{theorem}\label{thm:triangle-numbers} The maximum density of a packing of  $n$ equal disks $\delta_{\Delta}$ in a triangular torus is achieved when and only when $n = n_1^2 + n_1n_2 + n_2^2$ for $n_1, n_2$ integers.
\end{theorem}
Figure \ref{fig:triangular-packing} shows this for $n_1=2$, and $n_2=1$ and $n=7$.

We call any integer of the form (\ref{eqn:triangle-number}) a \emph{triangle lattice number}, which are also called Loeschian numbers in \cite{Conway-Sloane}, page 111, Sloane sequence $A003136$.  The first few numbers of this form are 0, 1, 3, 4, 7, 9, 12, 13, 16, 19, 21, 25, 27, 28, 31, 36, 37, 39, 43, 48, 49.  It is easy to see that the product of triangle lattice numbers is a lattice triangle number, and it is shown in \cite{Triangle-numbers} that $n>0$ is a triangle lattice number if and only if the non-square prime factors are congruent to $1$ modulo $6$, or it is the prime $3$. 
 
\subsection{The conjectures}\label{sect:conjectures}
 
One of our motivations is the following conjecture of L. Fejes Toth in \cite{solid-original,solid-conjecture}.  A packing of circles in the plane is said to be \emph{solid} if no finite subset of the circles can be rearranged so as to form, along with the rest of the circles, a packing not congruent to the original one.

\begin{conj}\label{conj:Fejes-Toth} (L. Fejes Toth) The triangle packing in the plane, minus one packing disk, is solid.
\end{conj}

The definition of a solid packing refers to unlabeled packings, so permuting the packing elements is not considered a distinct packing.  The following is a related conjecture, but in terms of finite packings in the triangular torus is what we call the Density Gap Conjecture.

\begin{conj}\label{conj:solid-finite} Let $n$ be a positive integer so that $n+1$ is a triangular lattice number but $n$ is not a triangular lattice number.  Then the most dense packing of $n$ equal disks in the triangular torus is $\delta = \frac{n}{n+1}\delta_{\Delta}=\frac{n}{n+1}\frac{\pi}{\sqrt{12}}$, and is achieved only by the triangular packing with one disk removed.
\end{conj}

The point is that Conjecture \ref{conj:solid-finite} implies Conjecture \ref{conj:Fejes-Toth} since a finite subpacking of a counterexample to the Conjecture \ref{conj:Fejes-Toth} can be placed in a triangular torus with a sufficiently large number of disks.  We do not provide a proof of Conjecture \ref{conj:solid-finite}, but we do provide some evidence for it in the following sections.  

Following the definitions in \cite{Donev1}, we say that a packing is \emph{collectively jammed} if the only continuous motion of the packing disks, fixing the lattice defining the torus and fixing the packing radius, is a continuous translation.  The following is another conjecture closely related to Conjecture \ref{conj:solid-finite}.

\begin{conj}\label{conj:local} If a packing of $n$ equal disks in a triangular torus is collectively jammed and is not a triangular packing, then its density is at most $\delta= \frac{n}{n+1}\delta_{\Delta}=\frac{n}{n+1}\frac{\pi}{\sqrt{12}}$, and the maximal density  is achieved only by the triangular packing with one disk removed. 
\end{conj}

Note that Conjecture \ref{conj:local} concerns any collectively jammed packing of disks, whether or not there is a triangular packing of $n+1$ disks available.

\subsection{Results for small numbers of packing disks}

\begin{figure}[here]
    \begin{center}
        \includegraphics[width=.7\textwidth]{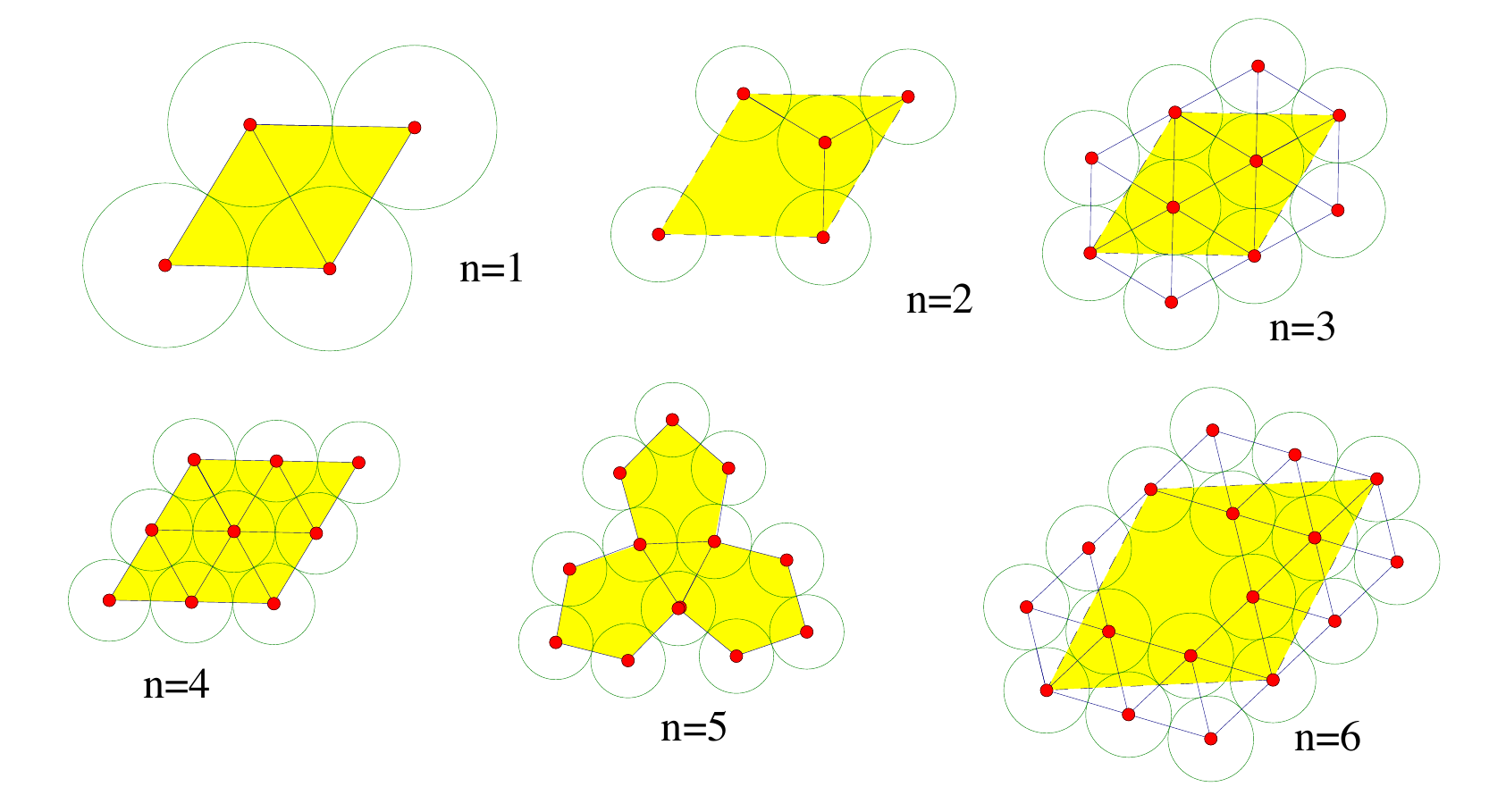}%
        \end{center}
    \caption{These are the locally maximally dense and maximally dense packings of equal disks in a triangular torus for $n$ from $1$ to $6$.  A fundamental region is shaded in each case.}
    \label{fig:best-packings}
    \end{figure}

In \cite{Will-students} the locally maximally dense packings of $n$ equal disks in a triangular torus were investigated, and all the locally maximally dense packings were found for $1 \le n \le 6$.  Figure \ref{fig:best-packings} shows those packings.  For each of those values of $n$, there is only one locally maximally dense packing, and thus it is the maximally dense packing. 
For all of these packings, except $n=5$ these are the packings expected from Conjecture \ref{conj:solid-finite}, namely, either a triangular lattice packing or such a packing minus one packing disk.  For five disks, a fundamental region is shown that has three-fold dihedral symmetry, and thus the whole packing has that symmetry.  Similar results were found in \cite{Will-square} for a square torus.

\subsection{Rigidity theory}
 
 For any packing the \emph{graph of the packing} is obtained by connecting each pair of centers of the packing disks that are tangent to each other with an edge.  Indeed the (labeled) disk centers form a configuration of points $\p=(\p_1, \p_2, \dots)$ in the plane or quotient torus, and we consider each edge $\{j,k\}$ as a strut in the sense of a tensegrity framework as in \cite{Connelly-Danzer-packing, Connelly-PackingI, Connelly-PackingII}.  An infinitesimal flex of framework $G(\p)$ is a sequence of vectors $\p' = (\p_1', \p_2', \dots)$ where each $\p_j'$ is associated to $\p_j$ such that for each strut $\{j,k\}$
\begin{equation}\label{eqn:strut-constraint}
(\p_j- \p_k) \cdot (\p_j' - \p'_k) \ge 0.
\end{equation}
Although a periodic configuration $\p$ is infinite, for a fixed lattice $\Lambda$, we consider each equivalence class $\p_j+ \Lambda$ as a single point in the quotient torus, and most of our calculations are done in the Euclidean plane.  The only continuous isometries of the quotient torus are translations, so the only trivial infinitesimal flexes of our torus tensegrities are when each $\p_j'$ are the same constant.  With this in mind, we say that a strut tensegrity (or more generally any tensegrity) is \emph{infinitesimally rigid} if the only infinitesimal flexes are trivial.  A strut tensegrity, in a torus, is \emph{rigid} if the only continuous motions of the vertices not decreasing the lengths of edges are translations.  The following is a useful result that can be found in \cite{Connelly-PackingI, Connelly-PackingII}:

\begin{theorem}\label{thm:strut-motion}Any strut tensegrity in a fixed torus is rigid if and only if it is infinitesimally rigid.
\end{theorem}

That infinitesimal rigidity implies rigidity is a general fact.  The converse implication is because all the edges are struts.  

If the infinitesimal flexes $\p'$ that satisfy (\ref{eqn:strut-constraint}) are equalities, we say that $\p'$ is an infinitesimal flex of the \emph{bar framework}.  Again if the only infinitesimal flexes of the bar framework are trivial, we say the bar framework is infinitesimally rigid.  

Any assignment of a scalar $\omega_{jk}=\omega_{kj}$ to each edge $\{j,k\}$ of the packing graph is called a \emph{stress} (where $\omega_{jk}=0$ when $\{j,k\}$ is not an edge in the graph) and it is called an \emph{equilibrium stress} or equivalently a \emph{self-stress}, if for each vertex $\p_j$
\begin{equation*}
\sum_k \omega_{jk}(\p_k - \p_j)= 0,
\end{equation*}

The following result of B. Roth and W. Whiteley \cite{Roth-Whiteley} is very useful in deciding when a strut tensegrity is infinitesimally rigid, because often solving the equality constraints of the bar conditions are somewhat easier than the direct inequality constraints of (\ref{eqn:strut-constraint}).

\begin{theorem}\label{thm:Roth-Whiteley}A strut tensegrity is infinitesimally rigid if and only if the underlying bar framework is infinitesimally rigid and there is an equilibrium stress that is negative on all the struts.
\end{theorem}

So the lack of an equilibrium stress, all of the same sign, or an infinitesimal motion of the bar framework, is enough to show that the disks can be moved keeping the packing property. 

With these results in mind we can model the local behavior of a packing of disks with a strut tensegrity on the edges connecting the centers of tangent packing disks.  If the underlying graph of a packing is rigid (and therefore infinitesimally rigid) as a strut graph in a fixed torus, then the packing is collectively jammed as defined in Subsection \ref{sect:conjectures}.

 A consequence of Theorem \ref{thm:strut-motion} and Theorem \ref{thm:Roth-Whiteley} is that if a packing of disks is collectively jammed, then $m \ge 2n - 2 +1 = 2n -1$, where $n$ is the number of packing disks, and $m$ is the number of contacts, i.e. the number of edges in the contact graph.

We say that a packing is \emph{locally maximally dense} if there is an epsilon $\epsilon > 0$ such that any $\epsilon$  perturbation of the packing to another packing (with possibly different common radius) has the same or less density.  Another consequence of the proof of Theorem \ref{thm:strut-motion} is that any locally maximally dense packing must have a subpacking that is rigid and thus infinitesimally rigid.  We call such a rigid subset of the packing disks a \emph{(rigid) spine} of the original packing.  In other words, the spine is collectively jammed.  In many cases the spine consists of the whole packing, or there are isolated packing disks not in the spine.  A packing disk, not in the spine is called a \emph{rattler}.  For example, in \cite{Donev1}  there are examples of locally maximally dense circle packings in a square torus, or a nearly square torus, that has rattlers.

The following is our main result, which is a special case of Conjecture \ref{conj:solid-finite}.  The next sections will provide a proof. 
\begin{theorem}\label{thm:main} Suppose that the graph of a collectively jammed packing of $n$ disks in a triangular torus consists of not all triangles.  Assume that either all the faces are rhombi, or the rhombi form a single strip, where the rest of the faces are triangles.  Then its density is $\delta < \frac{n}{n+1}\delta_{\Delta} =  \frac{n}{n+1}\frac{\pi}{\sqrt{12}}$.
\end{theorem}  

\section{Triangle and square tilings}\label{sect:TS-tilings}
 
\begin{figure}[here]
    \begin{center}
        \includegraphics[width=.7\textwidth]{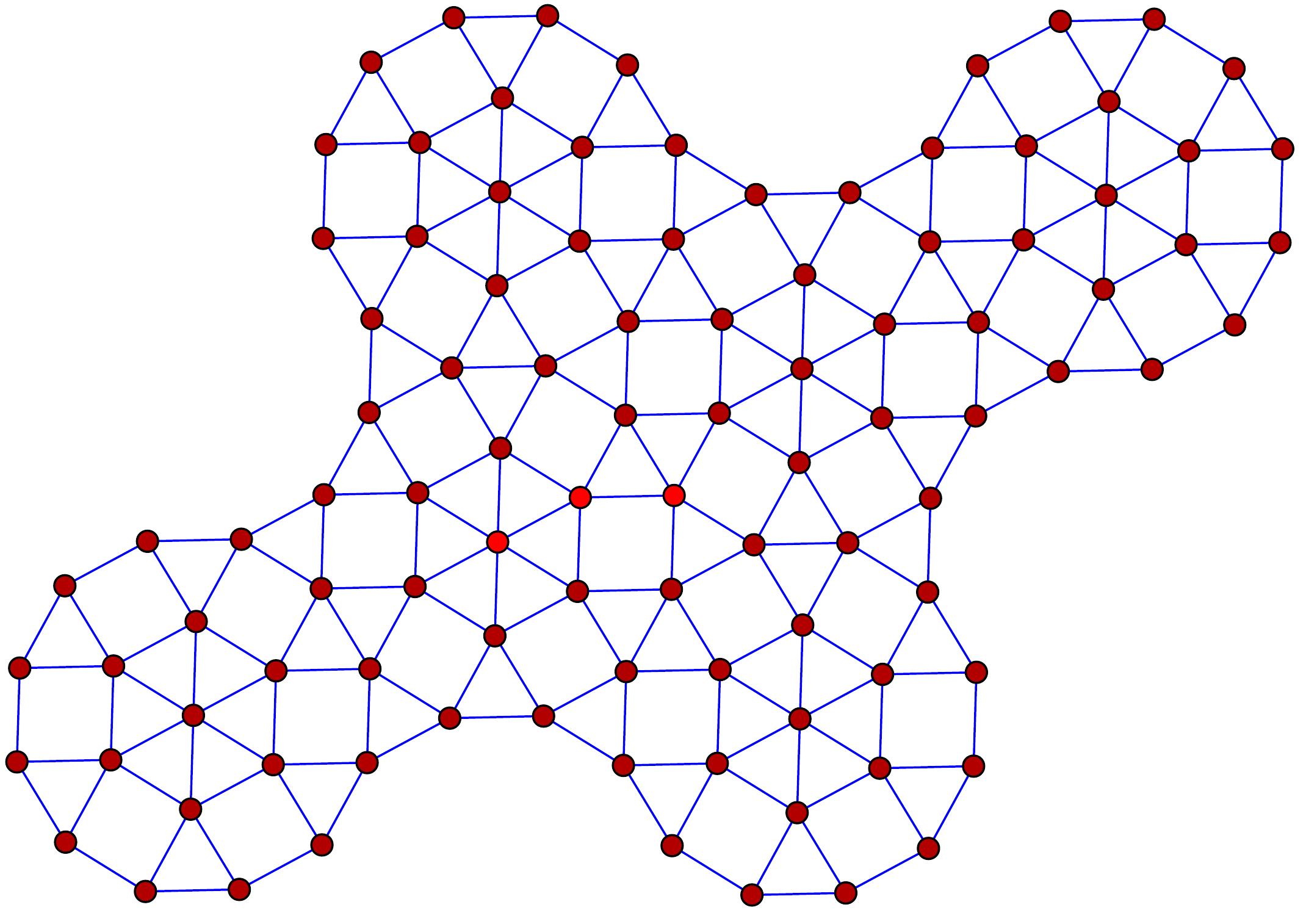}%
        \end{center}
    \caption{This shows a portion of a square and triangle tiling of the plane.}
    \label{fig:triangle-square-final-1}
    \end{figure}

In this section we provide information about the rigidity and flexibility properties of triangle and square tilings of the plane.  An example of a square and equilateral triangle tiling is shown in Figure \ref{fig:triangle-square-final-1}.

Given an equilateral triangle and square tiling of the plane, a circle packing is obtained by centering circles of radius one-half of the edge length at each of the vertices.   In fact, any graph corresponding to an edge-to-edge tiling of triangles and rhombi, whose edge lengths are the same, and such that the internal angles are all greater than or equal to $2\pi/6$ will correspond to a circle packing. 

Let $g = e^{2\pi i/12}$, a twelfth root of unity in the plane $\C$.  $g$ is one of the four roots of the cyclotomic polynomial $z^4 -z^2 +1$, which can be regarded as a defining relation for $g$, and all $12$ powers of $g$ are an integral linear combination of $1, g, g^2, g^3$.  So $\Lambda = \Lambda(1, g, g^2, g^3)$ is a four-dimensional lattice, and any point $z$ in $\Lambda$ is a unique integral linear  combination $n_1 + n_2 g+n_3 g^2 +n_4 g^3=z$, since the cyclotomic polynomial is irreducible over the rational field.  (Alternatively, this follows since $1$ and $\sqrt{3}/2$ are independent over the rationals, and they form the real and imaginary components of the four powers of $g$.)  

Position the origin at one of the vertices of the tiling and align one edge of graph with one of the twelve generators of $\Lambda$.  One sees that each oriented edge in the tiling is one of the powers of $g$ and each vertex of the tiling is in $\Lambda$.  Also $g^2$ is a sixth root of unity with cyclotomic polynomial $z^2 -z +1$, so $\Lambda_1=\Lambda(1, g^2)$ and $\Lambda_2=\Lambda(g, g^3)$ are both sublattices of $\Lambda$.  Each directed edge (directed in either direction) of each triangle in the triangle and square tiling belong to exactly one of $\Lambda_1$ or $\Lambda_2$.  We say that two polygons in the plane, with disjoint interiors, are \emph{edge-adjacent} or \emph{edge connected} if they have an edge in common.  We say they are \emph{vertex-adjacent} if they have a vertex in common, but are not edge-adjacent.  When two triangles are edge-adjacent, their edges belong to the same sublattice, and when they are vertex-adjacent to each other and edge-adjacent to the same square, their edges are in opposite lattices.  

\begin{figure}[here]
    \begin{center}
        \includegraphics[width=.5\textwidth]{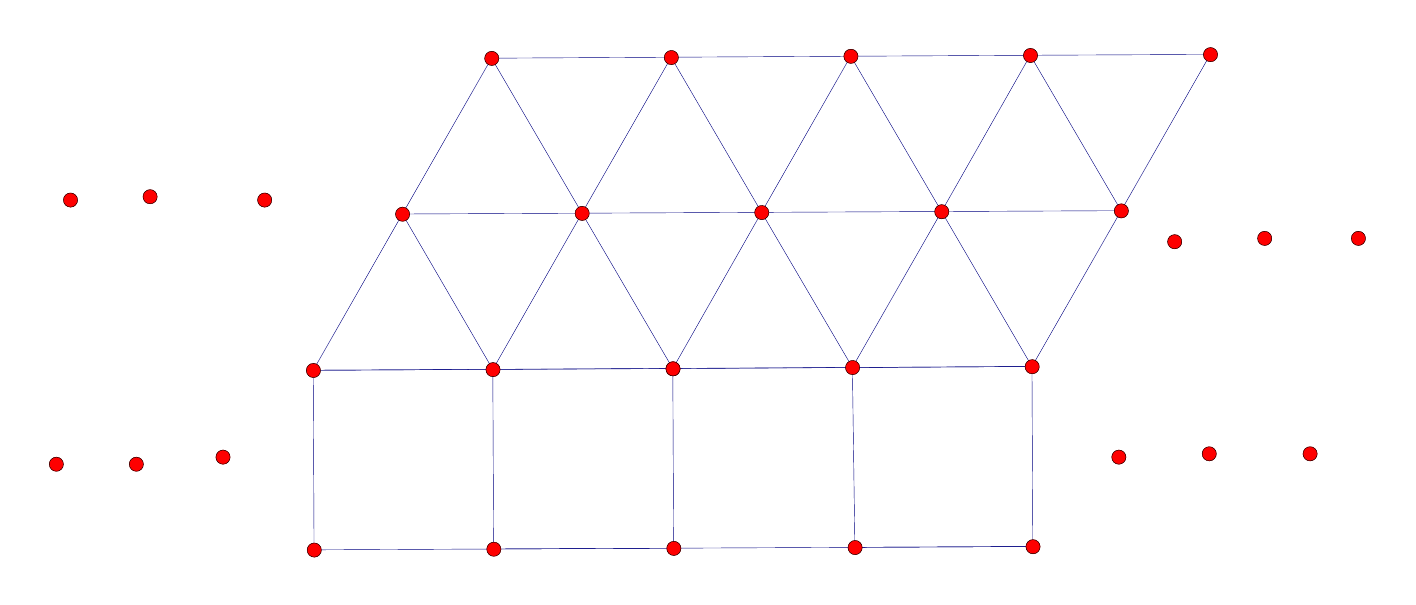}%
        \end{center}
    \caption{This shows a strip in a square-and-triangle tiling.}
    \label{fig:strips}
    \end{figure}

\begin{figure}[here]
    \begin{center}
        \includegraphics[width=.4\textwidth]{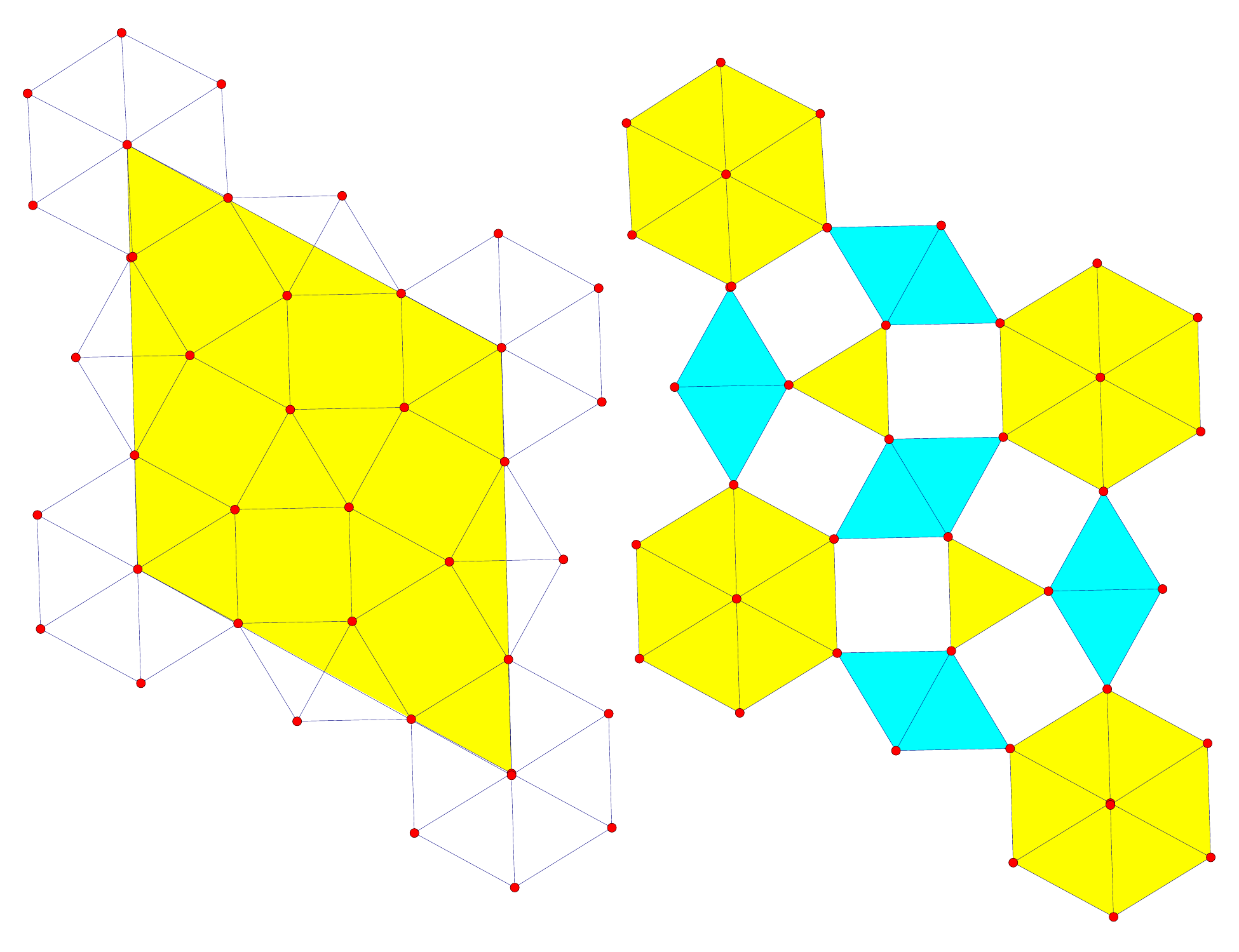}%
        \end{center}
    \caption{This shows a periodic triangle and square tiling with a period rhombus shown and the rigid groups of triangles shown.  The triangles with edges in the same lattice have the same shade.  The $3$-fold covering of this tiling is the one given in Figure \ref{fig:triangle-square-final-1}}
    \label{fig:triangle-square-flex}
    \end{figure}

\begin{figure}[here]
    \begin{center}
        \includegraphics[width=.9\textwidth]{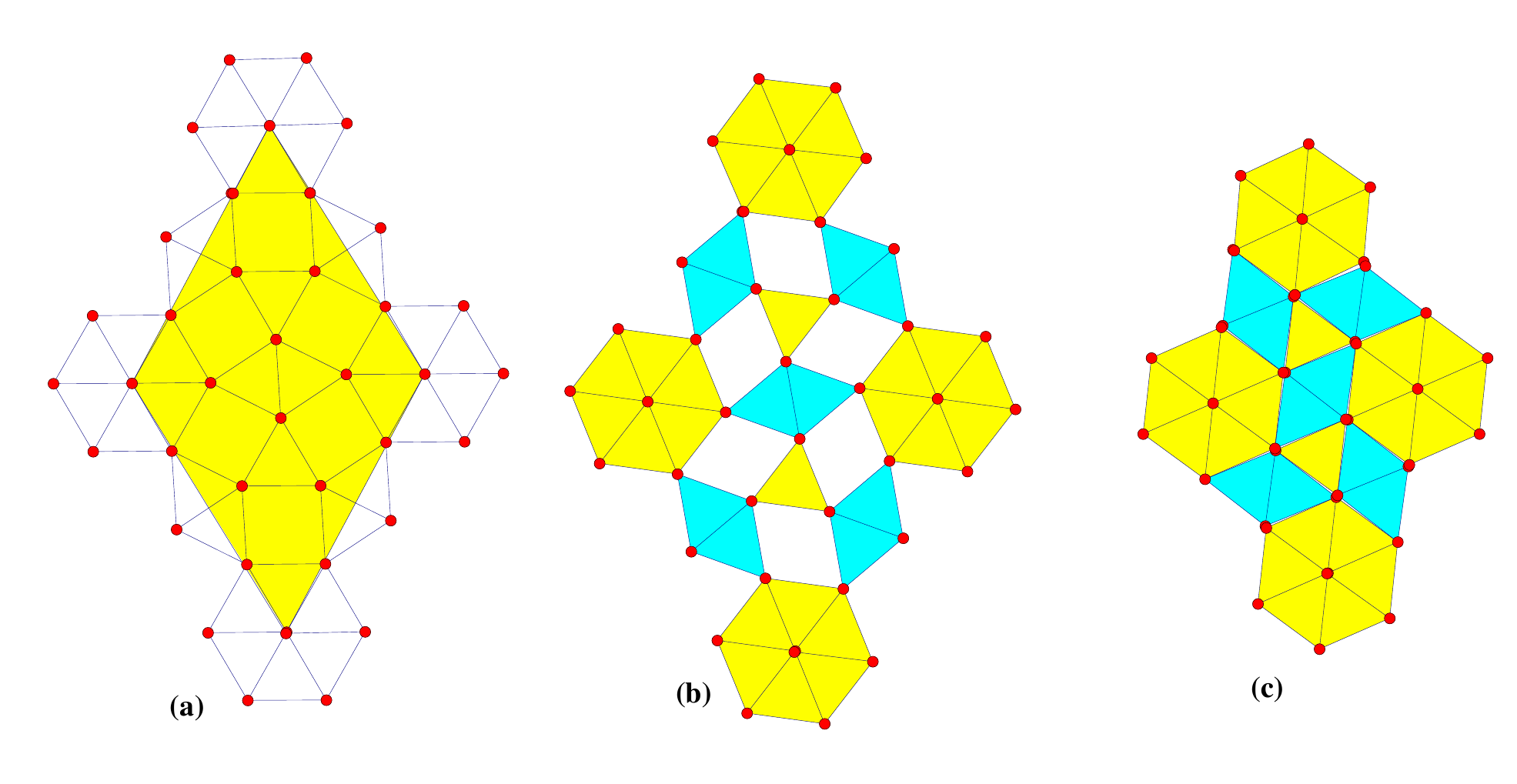}%
        \end{center}
    \caption{This shows the flexed framework of Figure \ref{fig:triangle-square-flex} where each square is flexed until it is the union of two triangles.  When the extra edge is filled in, the resulting tiling is a triangle lattice tiling.  There are $13$  vertices in the original tiling.  When it is flexed even further in Figure (c) to collapse the squares entirely, another tiling of triangles occurs, in this case with $7$ vertices.}
    \label{fig:triangle-square-flexed}
    \end{figure}

Consider an infinite edge-to-edge tiling of the plane by triangles and squares such that there are no infinite strips of squares or triangles, as in Figure \ref{fig:strips}.
\begin{theorem}\label{thm:main}The bar graph consisting of the edges of a square and equilateral  triangle tiling of the plane, with no infinite strips, flexes with one degree of freedom by rotating the edges of the lattice $\Lambda_1$ relative to the lattice $\Lambda_2$.  If the tiling is periodic with respect to a triangular lattice, it remains periodic with respect to a contracted triangular lattice.  If the tiling is periodic with respect to a square lattice it remains periodic, and the period lattice generators remain perpendicular, but they do not remain the same length.  
\end{theorem}
\begin{proof}
If an edge of the square and triangle tiling is in $\Lambda_1$, it must be an even power of $g$, and if it is in $\Lambda_2$, it is an odd power of $g$.  This is indicated by the shading in Figure \ref{fig:triangle-square-flex} and Figure \ref{fig:triangle-square-flexed}.  Note that both $\Lambda_1$ and $\Lambda_2$ are two-dimensional discrete lattices, so they can be rotated by $g_{\theta}=e^{i\theta}$ and the same integral linear combination of the generators of $\Lambda$ can be used for $\Lambda_1$ and $\Lambda_2$.  If there are no infinite strips, there must be edges of the tiling that are part of the generators from both $\Lambda_1$ and $\Lambda_2$.  So, explicitly, if $z = n_1 + n_2 g + n_3 g^2 + n_4 g^3$  is in $\Lambda$, then the flexed $z$ is
\begin{equation}\label{formula:turned-vector}
z(\theta) = n_1 + n_2 gg_{\theta} + n_3 g^2 + n_4 g^3g_{\theta}= \lambda_1 + \lambda_2   g_{\theta},
\end{equation}
where $ \lambda_1= n_1 + n_3 g^2$ is in $\Lambda_1$, and $\lambda_2 = n_2 g + n_4 g^3$ is in $\Lambda_2$.

The next question is to determine what the lattice of the flexed configuration is.  Let $\lambda_1 + \lambda_2$ represent one of the generators of a period parallelogram for a triangular lattice tiling of triangles and square, where $\lambda_1 \in \Lambda_1$ and $\lambda_2 \in \Lambda_2$.  Then 
$g^2(\lambda_1 + \lambda_2)$, the rotated vector by $2\pi/6$, is another generator.  Then in the flexed configuration these generators are flexed to $\lambda_1 + \lambda_2 g_{\theta}$, and $g^2(\lambda_1 + \lambda_2 g_{\theta})$, respectively, since $g^2\lambda_1 \in \Lambda_1$ and $g^2\lambda_2 \in \Lambda_2$, and clearly the length of the these generators are  $|\lambda_1 + \lambda_2 g_{\theta}| = |g^2(\lambda_1 + \lambda_2 g_{\theta})|$.   But for the same pair of generators $g^4(\lambda_1 + \lambda_2)$ is the difference vector between them, and by a similar calculation it also has the same length as the other generators.  Thus the lattice is a contracted triangular lattice, as desired.

For the square lattice, there is an important difference in how the generators are flexed.  If $\lambda_1 + \lambda_2$ is a generator as before, then $g^3(\lambda_1 + \lambda_2)$ is the other generator of the same length at right angles.  But now $g^3\lambda_1 \in \Lambda_2$, and $g^3\lambda_2 \in \Lambda_1$, so the flexed generators are $g^3(\lambda_1 g_{\theta} + \lambda_2)$ and $g^3(\lambda_1 + \lambda_2 g_{\theta})$, respectively.  But nevertheless, as long as there are no strips, the infinite framework can still be flexed as before, and when $g_{\theta}=g_{2\pi/6}$ all the squares have been flexed to $\pi/6$ rhombi, and the vertices are part of a triangular lattice.  If the period parallelogram remains a square, it means that square lattice is a sublattice of the triangular lattice.  But that is not possible, since if $\lambda \in \Lambda_1$, then its rotation by $2\pi/4$, which is $g^3\lambda$, would be in $\Lambda_2$.  But these two lattices are disjoint except at $0$.  $\Box$
\end{proof}


Note that a consequence of Theorem \ref{thm:main} is that that any packing whose graph consists of triangles and squares, with no strips in a triangular torus, is not rigid, i.e. not collectively jammed.   

\section{Triangle and rhombus tilings}\label{sect:triangle-rhombus}

We next consider the case when the quadrilaterals are not necessarily  squares, but rhombi.  For the case when the tiling is the graph of a packing of equal circles, we assume that each internal angle of each rhombus is strictly between $2\pi/6=\pi/3$ and $\frac{2}{6}2\pi=2\pi/3$.   If the angle is either of those extreme angles, we assume that the rhombus is subdivided into two equilateral triangles.  If the graph has that property, we say that it is \emph{completed}.

\begin{theorem}\label{thm:rhombi-triangle}Suppose a (edge-to-edge) tiling of the plane into (equilateral) triangles, rhombi, such that each vertex that is adjacent to a triangle has degree $5$ or $6$, and there are no infinite strips.  Then all the rhombi are congruent and the tiling can be flexed to a square and triangle tiling as in Theorem \ref{thm:main}.  
\end{theorem}
\begin{proof}  The edge-connected triangle components form convex polygons with between $3$ and $6$ vertices, with internal angles other than $\pi$, and such that the internal angle at each other vertex is either $2\pi/3$ or $2\pi/6$.  There are no non-convex vertices because of the completed property.  Figure \ref{fig:triangle-polygon} shows examples of these polygons for all the cases from $3$ to $6$.
\begin{figure}[here]
    \begin{center}
        \includegraphics[width=.8\textwidth]{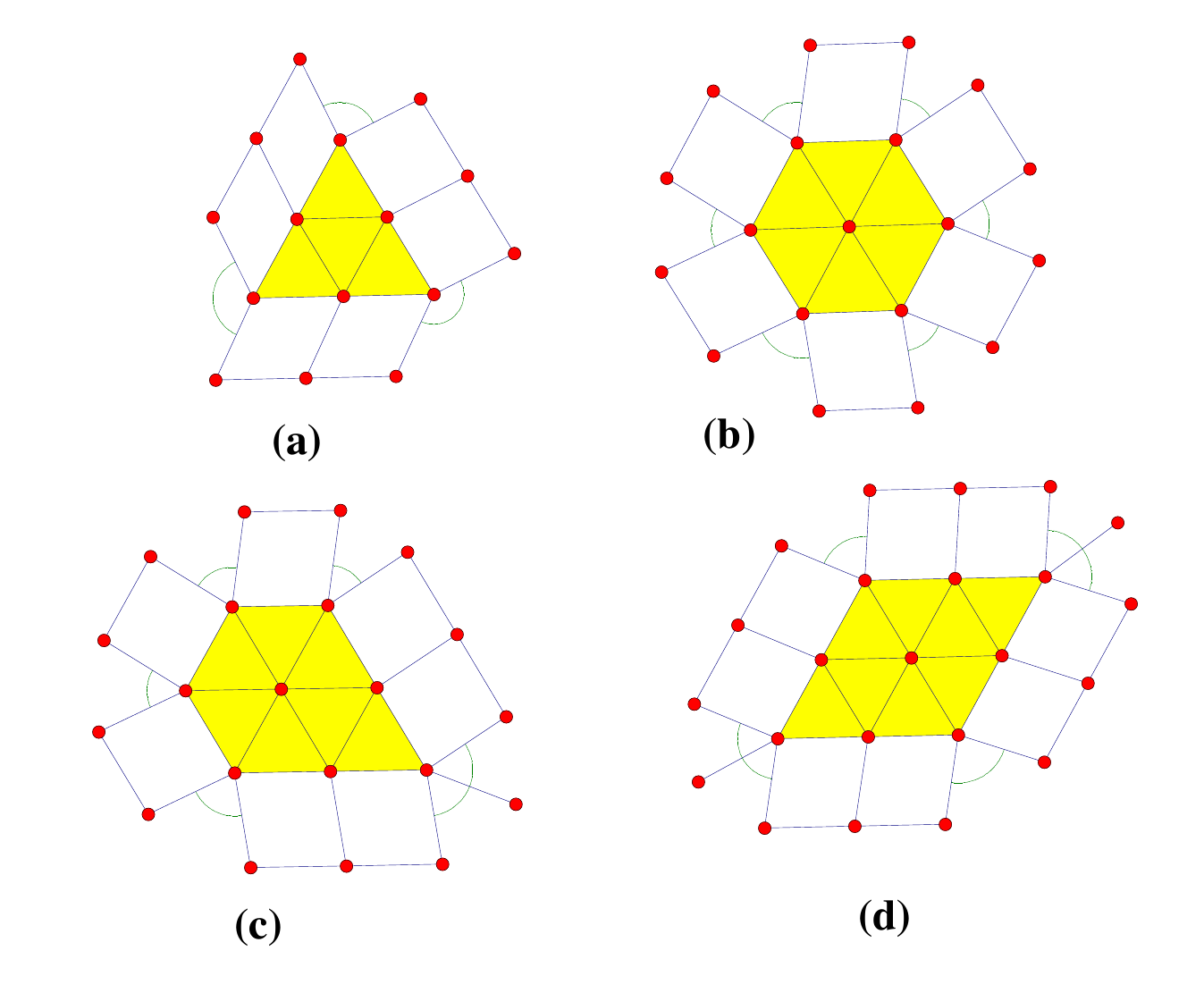}%
        \end{center}
    \caption{The shaded polygons are examples of components of triangles.  The edge-adjacent rhombi are shown as well.  The extra edges for cases (c) and (d), when the internal angles are $\pi/6$, are there because of the degree $5$ hypothesis.}
    \label{fig:triangle-polygon}
    \end{figure}
    
 On each edge  of each of these polygons, there is a rhombus or adjacent rhombi such that they form a parallelogram edge-adjacent to one edge of the triangle component.   The external angles of the triangle component are such that when you count the angles, not coming from the edge-adjacent rhombi, they add to $2\pi$.  In cases (b), (c), and (d) there are $6$ such angles, so the average angle is $2\pi/6=\pi/3$.  But since each angle must not be less than $\pi/3$, they are all $\pi/3$.  In case (a) there are three such angles, the average is $2\pi/3$ and each angle cannot be greater than $2\pi/3$.  So all the angles are $2\pi/3$.  Thus each triangle component is vertex-adjacent to another triangle component at each vertex (not including the vertices whose interior and exterior angles are $\pi$).  Furthermore,  as one proceeds around the boundary of each triangle component, the rhombi are all congruent, with the same internal angles (when they are not squares), appearing in the same order with respect to a clockwise order, say, along the boundary.   So the rhombus components are all parallelograms with those same angles, since there are no infinite components.  So each triangle component has its edge directions in one of two possible sets of directions, where each set consists of three directions $2\pi/3$ apart.  Thus one set can rotate relative to the other.  
 $\Box$\end{proof}
 
\begin{figure}[here]
    \begin{center}
        \includegraphics[width=.8\textwidth]{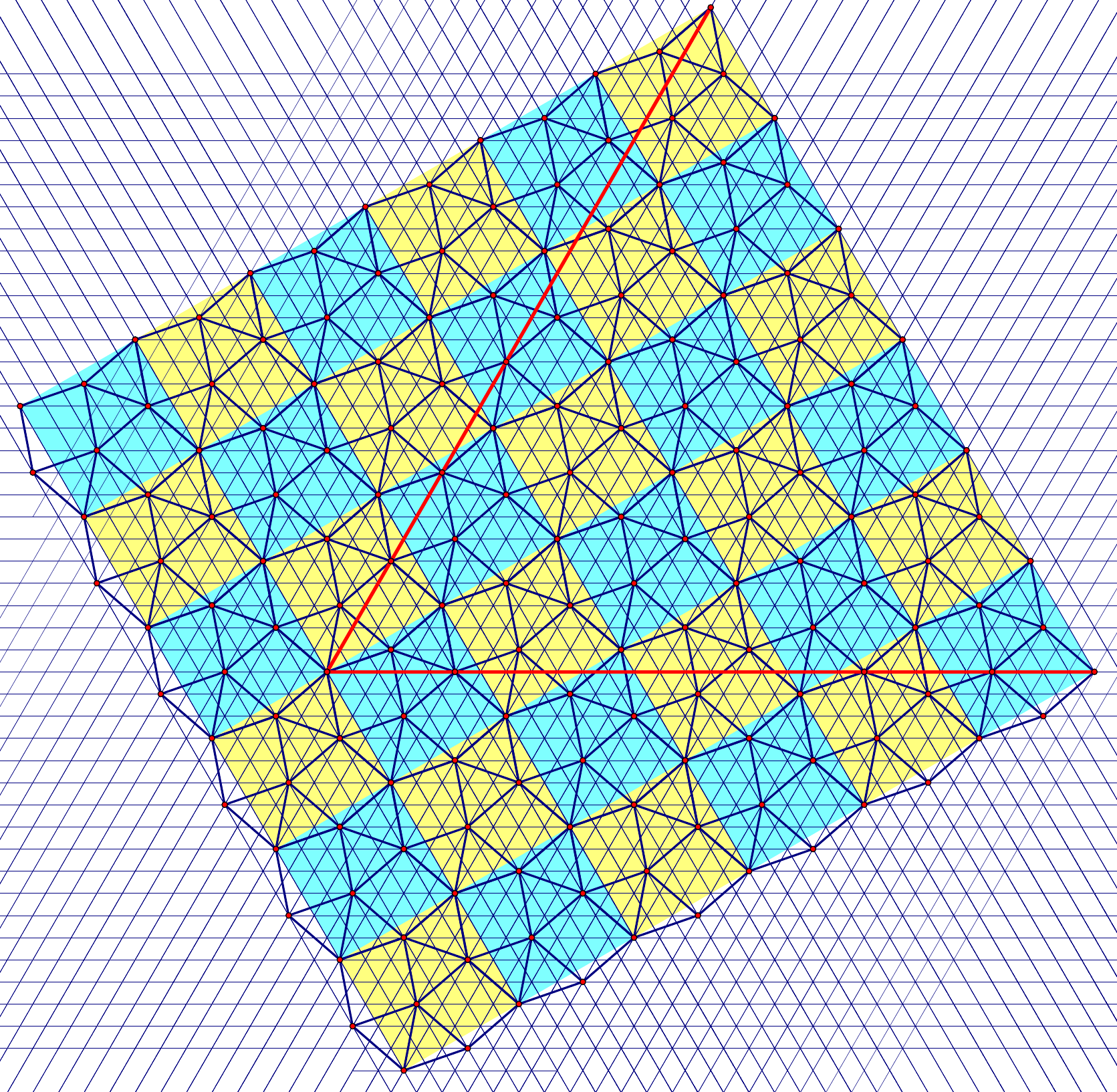}%
        \end{center}
    \caption{Two solid line segments here represent the generators of the $\Lambda_{\Delta}$ lattice.}
    \label{fig:Counterexample-1}
    \end{figure}

 \bigskip
\noindent{\bf Remark:} In a previous draft of this paper it was claimed that if the original tiling was periodic with respect to a triangular lattice, the the flex of Theorem \ref{thm:rhombi-triangle} would preserve this triangular periodicity.  This is not necessarily true.  Figure \ref{fig:Counterexample-1} shows an example.  The vertices of the packing graph lie on an underlying triangular lattice $\Lambda_{\epsilon}$ (which must be true for this sort of example) as shown in Figure \ref{fig:Counterexample-1}.  The lattice $\Lambda_{\epsilon}$ provides coordinates for the packing graph, and it allows one to verify that the configuration exists as shown.  The whole graph is a covering graph of a four-vertex graph on another lattice $\Lambda$ that is not a triangular lattice.  This graph is a triangle and rhombus tiling, but the rhombi are not squares, although the rhombus shape is close to a square.  The triangular sublattice $\Lambda_{\Delta}$ of $\Lambda$ is indicated in  Figure \ref{fig:Counterexample-1} 
and again the $\Lambda_{\epsilon}$ lattice provides coordinates to verify that $\Lambda_{\Delta}$ a triangular lattice.  But when the packing graph is flexed as in Theorem \ref{thm:rhombi-triangle}, it is easy to verify that the $\Lambda_{\Delta}$ lattice ceases to be triangular.  We thank Alex Smith and Jeff Shen for pointing out the flaw in the proof of our previously incorrect statement of Theorem \ref{thm:rhombi-triangle}.

    
\section{Rigid rhombus and triangle strip tilings}\label{sect:rhombus-tiling}

When all the faces of packing graph are rhombi, we provide an explicit description what they look like in Subsection \ref{subsect:rigid-rhombus-tiling}.  Then in Subsection \ref{subsect:rhombus-density} we calculate the maximum density of such packings for a triangular torus.

\subsection{A description of a collectively jammed rhombus packings}\label{subsect:rigid-rhombus-tiling}

Another case of tilings corresponding to collectively jammed packing graphs in a triangular torus is when all of the faces are rhombi.  We first determine the abstract graph of such a packing.  For any edge-to-edge tiling of a torus let $v_j$ be the number of vertices of degree $j$ for $j = 3, 4,  \dots$, and let $f_j$ be the number of faces with $j$ sides, for $j = 3, 4, \dots$.  Define the averages
\begin{equation*}\label{eqn:averages}
\bar{v} = \frac{1}{v}\sum_3^{\infty} j v_j, \,\,\, \bar{f} = \frac{1}{f}\sum_3^{\infty} j f_j,
\end{equation*}
where $v = \sum_3^{\infty} v_j$ is the total number of vertices, and $f = \sum_3^{\infty} f_j$ is the total number of faces.  (Notice that all the sums are finite since there are only a finite number of vertices and faces.  The bar notation, used here for averages, is not to be confused with complex conjugation.)  Let $e$ be the total number of edges in the graph.  Then counting the edge-vertex adjacencies and the edge-face adjacencies we get
\begin{eqnarray*}
v\bar{v} = 2e = f \bar{f}.
\end{eqnarray*}
Assuming that each face is simply connected,  and since the Euler characteristic of a torus is $0$, we get 
 \begin{eqnarray*}
0=v - e + f &=&{2e}/{\bar{v}} -e +  {2e}/{\bar{f}}, \,\,\,\,\,\text{and} \\  {1}/{\bar{v}} + {1}/{\bar{f}} &=& {1}/{2}.
\end{eqnarray*}
So when all the faces are quadrilaterals, $\bar{f}=4$, and thus $\bar{v}= 4$ as well.  But if there are any vertices of degree $3$ and the graph is a packing graph, the three angles at any such vertex must each be $2\pi/3$, implying that there is a triangle in the graph.  Thus $v_3=0$, and so all the vertices must be of degree $4$.  Summing up, we have shown:
\begin{prop} Any packing graph of equal circles in a torus, where each face is a rhombus, has each vertex of degree four.
\end{prop}
So far we have not assumed that the torus is a triangular torus, or even that the packing is jammed in any sense.  We next consider when the packing is collectively jammed.

\begin{figure}[here]
    \begin{center}
        \includegraphics[width=.4\textwidth]{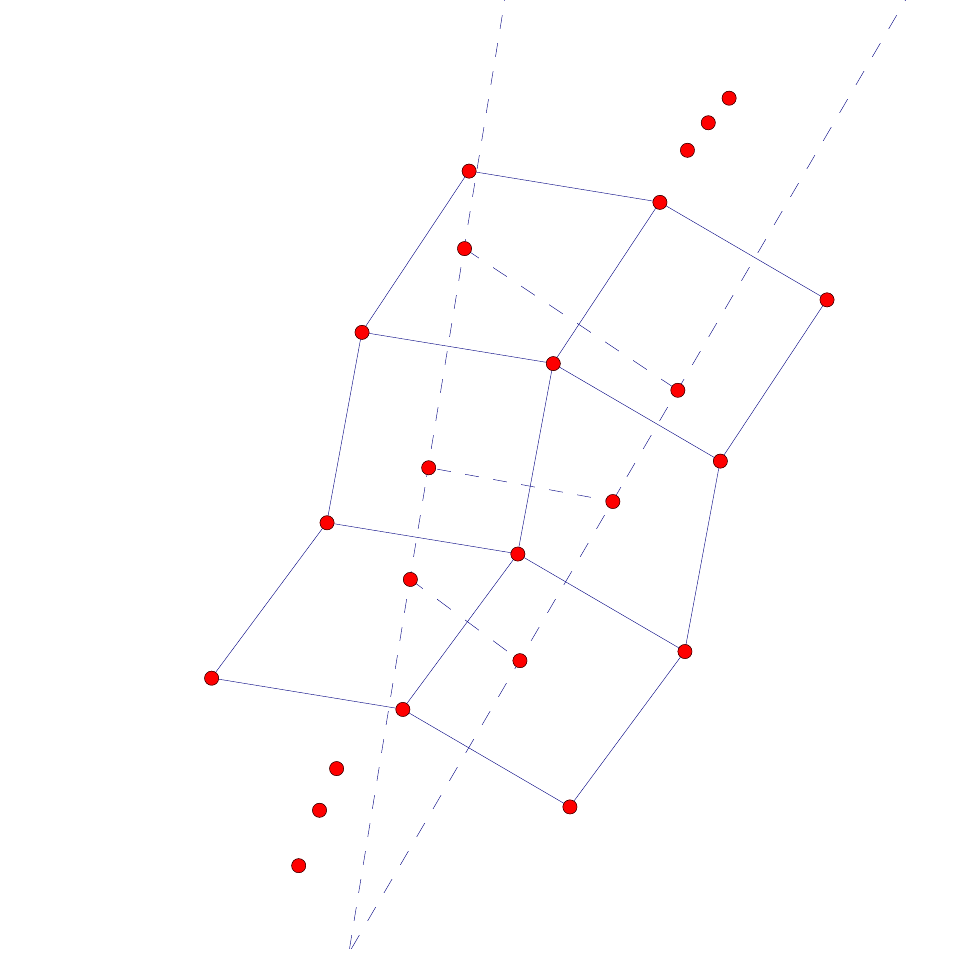}%
        \end{center}
    \caption{This shows a portion of an infinite rhombus tiling in the plane, where an attempt  is started to find a reciprocal graph is indicated in dashed lines, but cannot be continued.}
    \label{fig:parallel-strips}
    \end{figure}

If the graph of a collectively jammed packing consists of only rhombi, then in the lift to the plane, there are infinite strips of rhombi, where two opposite sides of each rhombus are in a fixed direction.  The opposite sides of each rhombus are adjacent to the next rhombus in the strip.  The other pair of sides are each adjacent to another rhombus in another strip of all different rhombi, where each rhombus in the second strip is edge-adjacent to a rhombus in the first strip.  Figure \ref{fig:parallel-strips} shows an example of such strips.  
 
By Theorem \ref{thm:strut-motion} and Theorem \ref{thm:Roth-Whiteley}, a necessary condition on the graph of a packing to be collectively jammed is that there is an equilibrium stress on the edges that is strictly negative on each edge.  One way of picturing this stress is to draw the \emph{reciprocal graph}, which is another graph whose vertices correspond to the faces of the original graph, whose faces correspond to the vertices of the original graph, and whose edges correspond to the edges of the original graph and are perpendicular to the original edges.  This reciprocal graph is indicated in Figure \ref{fig:parallel-strips} with dashed lines.  In the case when all the edges of the original graph are the same length, the stress on an edge corresponds to (minus) the length of the edge of the reciprocal.  See \cite{Connelly-PackingI, Connelly-PackingII} for a discussion of this with regard to packings.

Since the packing is periodic, the stress and the reciprocal extends to the infinite packing in the lift in the plane.  But each reciprocal polygon in the lift must be a convex polygon, which in our case is another quadrilateral, and the straight lines in the reciprocal, corresponding to the parallel sides of each adjacent strip, must not cross.  So this means that the lines are parallel.  Doing this for all non-crossing strips gives the following:

\begin{prop}\label{thm:parallel-rhombi}  If the graph of any collectively jammed circle packing in a torus consists of only rhombi, then the rhombi are all congruent and translates of a single rhombus.
\end{prop}

Proposition \ref{thm:parallel-rhombi} supplies the condition for the graph of the packing to have a proper stress that is necessary according to Theorem \ref{thm:Roth-Whiteley}, but we still need conditions for the graph to be infinitesimally rigid as a bar graph.

\begin{figure}[here]
    \begin{center}
        \includegraphics[width=.4\textwidth]{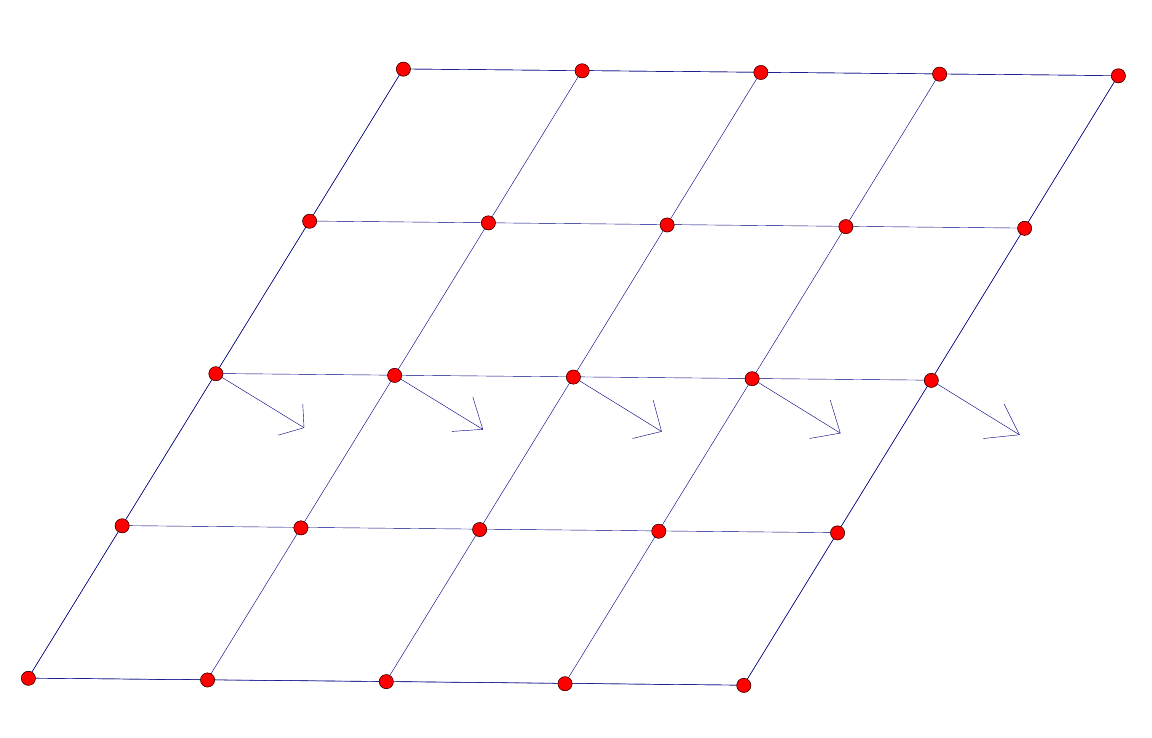}%
        \end{center}
    \caption{This shows an example of a periodic flex when the tour of centers is not all the vertices of the graph of a rhombus packing.}
    \label{fig:periodic-flex}
    \end{figure}

We now know that each edge of the rhombus tiling is in one of two directions.  As one proceeds from one vertex to the next on the torus in one of those directions the vertices are visited one after the other until one gets back to the starting vertex.  Suppose that $p_1, p_2, \dots, p_k$ is this \emph{tour} of vertices and it does not include all the vertices on the torus.  Then there is a non-zero infinitesimal flex $p'_1=p'_2= \dots= p'_k$ of the bar graph that is $0$ on all the other vertices not on this tour, where each $p'_j$ is perpendicular to the direction of the other edge of the original rhombus.  Figure \ref{fig:periodic-flex} shows an example of this.
It is easy to check that this is the only infinitesimal flex possible.  Hence we get the following:

\begin{theorem}\label{thm:tour} If the packing graph of a circle packing in a torus consists of only rhombi, the packing is collectively jammed if and only if all the rhombi are congruent with the sides in only two directions and the tour of vertices in each direction consists of all the vertices in the packing.
\end{theorem}

\begin{figure}[here]
    \begin{center}
        \includegraphics[width=.5\textwidth]{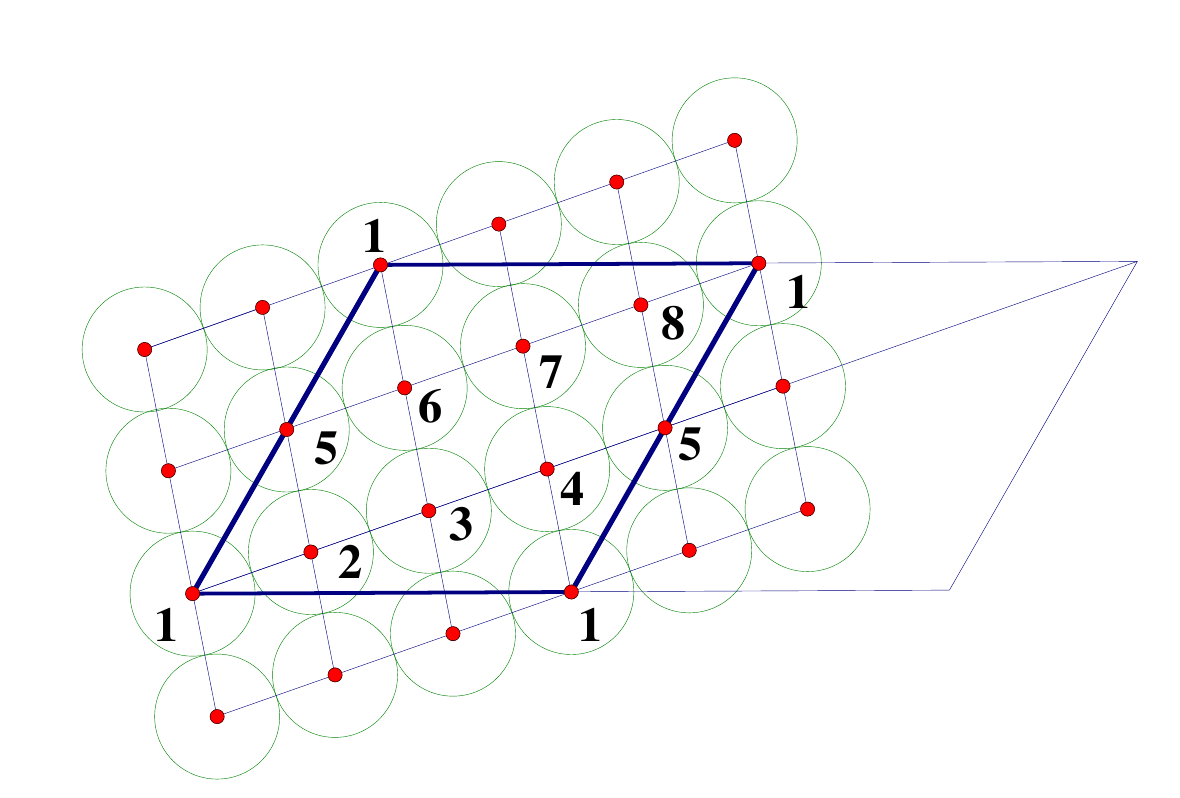}%
        \end{center}
    \caption{This shows an example of a collectively jammed packing on a triangular torus that, where one tour of vertices is indicated.  The other tour is $1, 4, 7, 2, 5, 8, 3, 6$.}
    \label{fig:rigid-rhombus}
    \end{figure}

Figure \ref{fig:rigid-rhombus} shows an example of the case when the packing graph consists of rhombi and each tour in each direction visits each vertex. 

\subsection{Strips of rhombi and triangles}\label{subsect:rhombus-triangles}
   
\begin{figure}[here]
    \begin{center}
        \includegraphics[width=.6\textwidth]{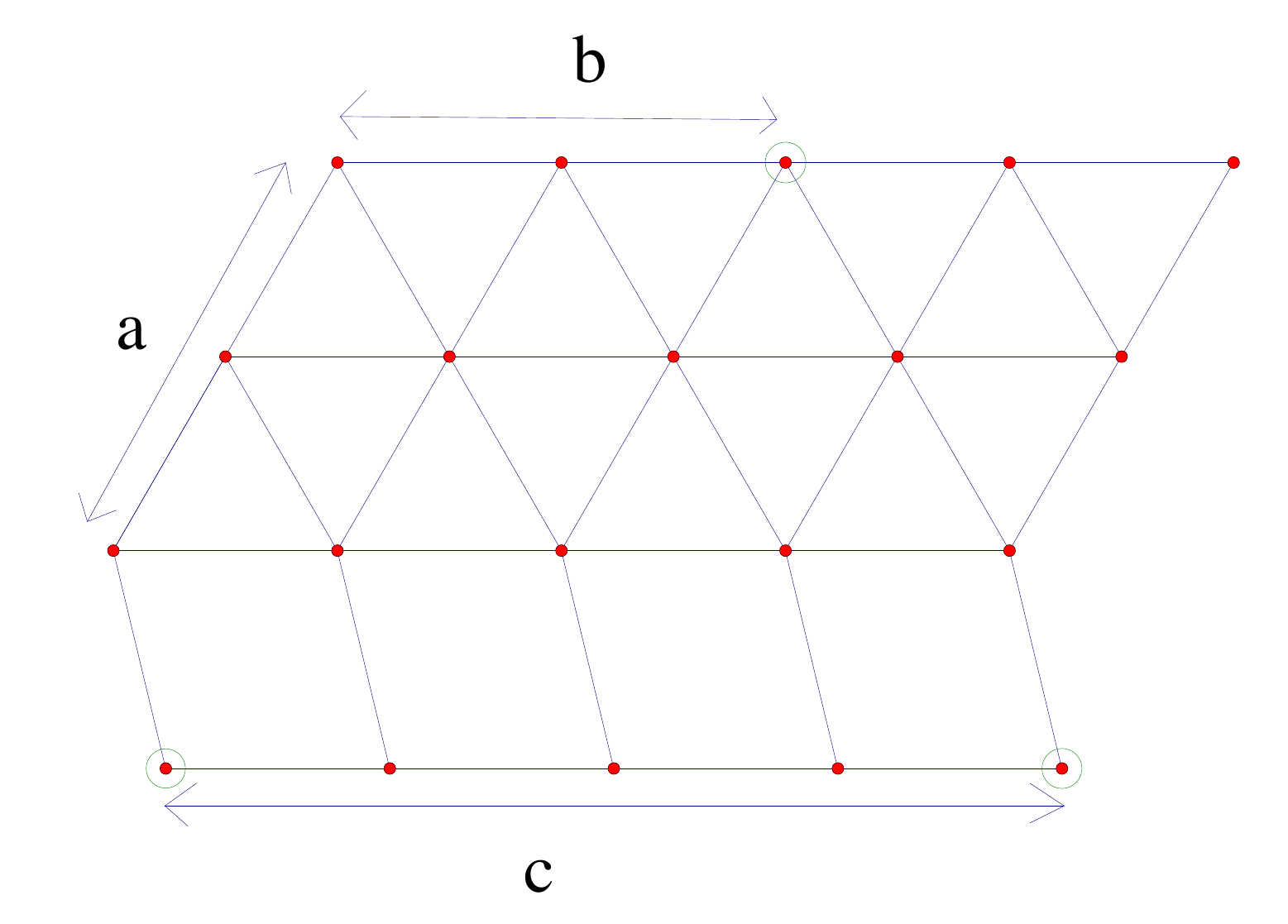}%
        \end{center}
    \caption{This shows a strip tiling of triangles and rhombi, where $c$ is the length of the strip of rhombi, $a$ is the width of strip of triangles, and $b$ represents the displacement of the identification of the bottom to the top.  Three of the lattice points are indicated.}
    \label{fig:strip-tiling}
    \end{figure}

 Another case of collectively jammed packings in a triangular torus is when there is a single strip of rhombi, and a corresponding strip of triangles of arbitrary thickness in the graph of the packing.  This includes, as a special case, the situation in Subsection \ref{subsect:rigid-rhombus-tiling} of all rhombi.  We call a strip of rhombi in a tiling of a \emph{torus strip} if each  rhombus is edge-adjacent with exactly two other rhombi and all lines containing the intersections between rhombi are parallel when lifted to the plane.  So these rhombi form a cyclic chain in the torus. In the case here, except for the case when all the faces are rhombi, we assume that there is only one such strip of rhombi, and all the other faces of the graph are (equilateral) triangles.  This is indicated in Figure \ref{fig:strip-tiling}.  In Section \ref{subsect:rigid-rhombus-tiling}, we assumed that all the faces were rhombi, and there were only two strips of rhombi, each one corresponding one of the two parallel classes of edges in the lift.  If there are two parallel strips of rhombi, it is easy to show that the graph would not be infinitesimally rigid and the packing would not be collectively jammed.
 
Place one of the lattice vectors, the origin, at lower left corner of one of the rhombi.  Then as you proceed in the direction of the strip, you must come back to the same point on the torus for the first time, and thus you must hit another lattice point in the lift to the plane.  Let $a$ be the number of layers of triangles, and $b$ the number of steps needed before you hit another lattice point on the upper layer, as indicated in Figure \ref{fig:strip-tiling}.  Thus the total number of vertices on the torus is $n= c(a+1)$.  The case when $a=0$ is the situation when all the faces are rhombi as in Subsection \ref{subsect:rigid-rhombus-tiling}.  

In order to do calculations of density, we orient the packing so that the strips all are aligned along the real axis, the edges all have unit length,  the rhombus edge direction, not along the real axis, is given by the unit length complex number $g$, where $\text{Im}(g) > 0$, and the internal angle at the origin for the rhombus is between $\pi/3$ and $2\pi/3$.  The edge directions of the triangles are $g_{\Delta}$ and $g_{\Delta}^2$ as before.
    
\subsection{Calculating the density of triangle-rhombus strip packings}\label{subsect:rhombus-density}
\begin{theorem} \label{thm:best-density}
The density $\delta$ of an equal circle packing whose graph is a triangle-rhombus strip tiling in a triangular torus is always such that $\frac{n}{n+1} \delta_{\Delta} > \delta$, where $n$ is the number of packing disks.
\end{theorem}

\begin{proof}
Since the area of the unit length rhombus is $\text{Im}(g)$, the total area of the torus is $(\text{Im}(g) + a\sqrt{3}/2)c$ and the density of the corresponding packing is
\begin{equation}\label{eqn:density}
\delta = \frac{\pi (1/2)^2 n}{(\text{Im}(g) + a\sqrt{3}/2)c} = \frac{\pi}{4}\frac{c(a+1)}{(\text{Im}(g) + a\sqrt{3}/2)c} = \frac{\pi}{4}\frac{a+1}{\text{Im}(g) + a\sqrt{3}/2}.
\end{equation}
The next task is to use the condition that the packing is in a triangular torus.  Let $v$ be a generating vector for a lattice that defines the triangular torus, so that lattice is $\Lambda(v, g_{\Delta}v)$.  Since $\Lambda(c, g + ag_{\Delta} +b)$ is the same lattice, we have
\begin{eqnarray*}
c &=& n_1 v + n_2 g_{\Delta}v = (n_1  + n_2 g_{\Delta})v\\
g + ag_{\Delta} +b &=& n_3 v + n_4 g_{\Delta}v = (n_3  + n_4 g_{\Delta})v,
\end{eqnarray*}
where $n_1, n_2, n_3, n_4$ are integers.  Then paying attention to the order of the lattice generators
\begin{equation}\label{eqn:det}
\det \begin{pmatrix}
n_1 & n_2 \\
n_3 & n_4
\end{pmatrix} = n_1 n_4 -n_2 n_3 = 1.
\end{equation}
Then we can solve for $v$ and simplify,
\begin{equation*}
v = \frac{c}{n_1  + n_2 g_{\Delta}} = \frac{c(n_1  + n_2 g_{\Delta}^{-1})}{n_1^2 + n_1n_2 +n_2^2}.
\end{equation*}
Similarly, we solve for $g$,
\begin{eqnarray}
g &=& -ag_{\Delta} -b + \frac{c(n_1  + n_2 g_{\Delta}^{-1})(n_3  + n_4 g_{\Delta})}{n_1^2 + n_1n_2 +n_2^2}\label{eqn:g-def} \\
   &=& -ag_{\Delta} -b + \left (\frac{c}{n_1^2 + n_1n_2 +n_2^2}\right )(n_1 n_3 +n_2 n_4 +n_1n_4g_{\Delta} + n_2 n_3 g_{\Delta}^{-1}).\nonumber 
\end{eqnarray}
We are primarily interested in the imaginary part of $g$, which is
\begin{eqnarray*}
\text{Im}(g) &=& \left[-a + \frac{c}{n_1^2 + n_1n_2 +n_2^2}(n_1n_4 -n_2n_3)\right] \frac{\sqrt{3}}{2}\\
&=& \left[-a + \frac{c}{n_1^2 + n_1n_2 +n_2^2}\right] \frac{\sqrt{3}}{2}
\end{eqnarray*}
by (\ref{eqn:det}).
Thus,
\begin{equation*}
\text{Im}(g) + a\frac{\sqrt{3}}{2} = \frac{c}{n_1^2 + n_1n_2 +n_2^2} \frac{\sqrt{3}}{2}.
\end{equation*}
Using (\ref{eqn:density}) we calculate the density of the corresponding packing,
\begin{eqnarray*}
\delta &=& \frac{\pi}{4} \frac{a+1}{\left(\frac{c}{n_1^2 + n_1n_2 +n_2^2}\right)\frac{\sqrt{3}}{2}} = \frac{\pi}{2\sqrt{3}} \frac{(a+1)(n_1^2 + n_1n_2 +n_2^2)}{c} \\
&=& \delta_{\Delta}\frac{(a+1)(n_1^2 + n_1n_2 +n_2^2)}{c}.
\end{eqnarray*}
Since the packings coming from strip tilings are not a triangular lattice packing $ \delta < \delta_{\Delta}$, and thus
\begin{equation}\label{eqn:gross-bound}
c > (a+1)(n_1^2 + n_1n_2 +n_2^2) \ge n_1^2 + n_1n_2 +n_2^2,
\end{equation}
and 
\begin{equation*}
\frac{ n_1^2 + n_1n_2 +n_2^2}{c} < 1.
\end{equation*}
Then we can improve the inequality (\ref{eqn:gross-bound}) to 
\begin{eqnarray*}
c &>& (a+1)(n_1^2 + n_1n_2 +n_2^2)+ \frac{ n_1^2 + n_1n_2 +n_2^2}{c}\\
c^2 &>& c(a+1)(n_1^2 + n_1n_2 +n_2^2) + n_1^2 + n_1n_2 +n_2^2\\
c^2 &>& [c(a+1) +1] (n_1^2 + n_1n_2 +n_2^2).
\end{eqnarray*}
Then going back to the total number of packing elements $n=c(a+1)$, we get 
\begin{equation*}
\frac{n}{n+1} = \frac{c(a+1)}{c(a+1) + 1} > \frac{(a+1)(n_1^2 + n_1n_2 +n_2^2)}{c}.
\end{equation*}
Thus 
\begin{equation*}
\frac{n}{n+1} \delta_{\Delta} > \delta_{\Delta}\frac{(a+1)(n_1^2 + n_1n_2 +n_2^2)}{c} = \delta.
\end{equation*}
\end{proof} $\Box$

Note that it is not used that such triangle-rhombus strip packings are collectively jammed, but they are.  It is easy to find a self-stress on the (strut) framework.  For the edges that are in the direction of $g_{\Delta}$ choose a positive stress $\omega_1$, and for the direction $g_{\Delta}^2$ choose $\omega_2$ and for the direction $g$ choose $\omega_3$ such that $\omega_1g_{\Delta} + \omega_2g_{\Delta}^2 +  \omega_3(-g) = 0$.  Then any constant positive stress on the horizontal members will do.  

When $a>0$, the connected component of triangles is infinitesimally rigid all by itself, since it is connected in the torus.  The case when $a=0$ was essentially handled in Subsection \ref{subsect:rigid-rhombus-tiling} and Equation (\ref{eqn:det}).

So Theorem \ref{thm:best-density} verifies the density gap conjecture Conecture \ref{conj:local} in the case of strips of rhombi and triangles.

In order believe that these calculations correspond to actual packings on a triangular torus, we present a specific example for all values of $a = 0, 1, \dots$.  For any non-negative integer $a$, let $b=4a + 5$, $c=(a+1)7+1=7a+8$, $n_1=2$, $n_2=1$, $n_3=1$, and $n_4=1$.  Then 
\begin{eqnarray*}
n_1n_4-n_2n_3 &=& 1\\
n_1^2 + n_1n_2 + n_2^2 &=& 7 \\
v &=& \frac{7a+8}{7}(2+g_{\Delta}^{-1})\\
g &=& \frac{1}{7} + \frac{8}{7}\frac{\sqrt{3}}{2}i.
\end{eqnarray*}
One can verify that Equation (\ref{eqn:g-def}) holds for any $a \ge 0$.  Thus this corresponds to collectively jammed strip packings.  Figure \ref{fig:rigid-rhombus} corresponds the case when $a=0$.

\section{Relation to prevoius work}\label{sect:previous-work}

\begin{figure}[here]
    \begin{center}
        \includegraphics[width=.8\textwidth]{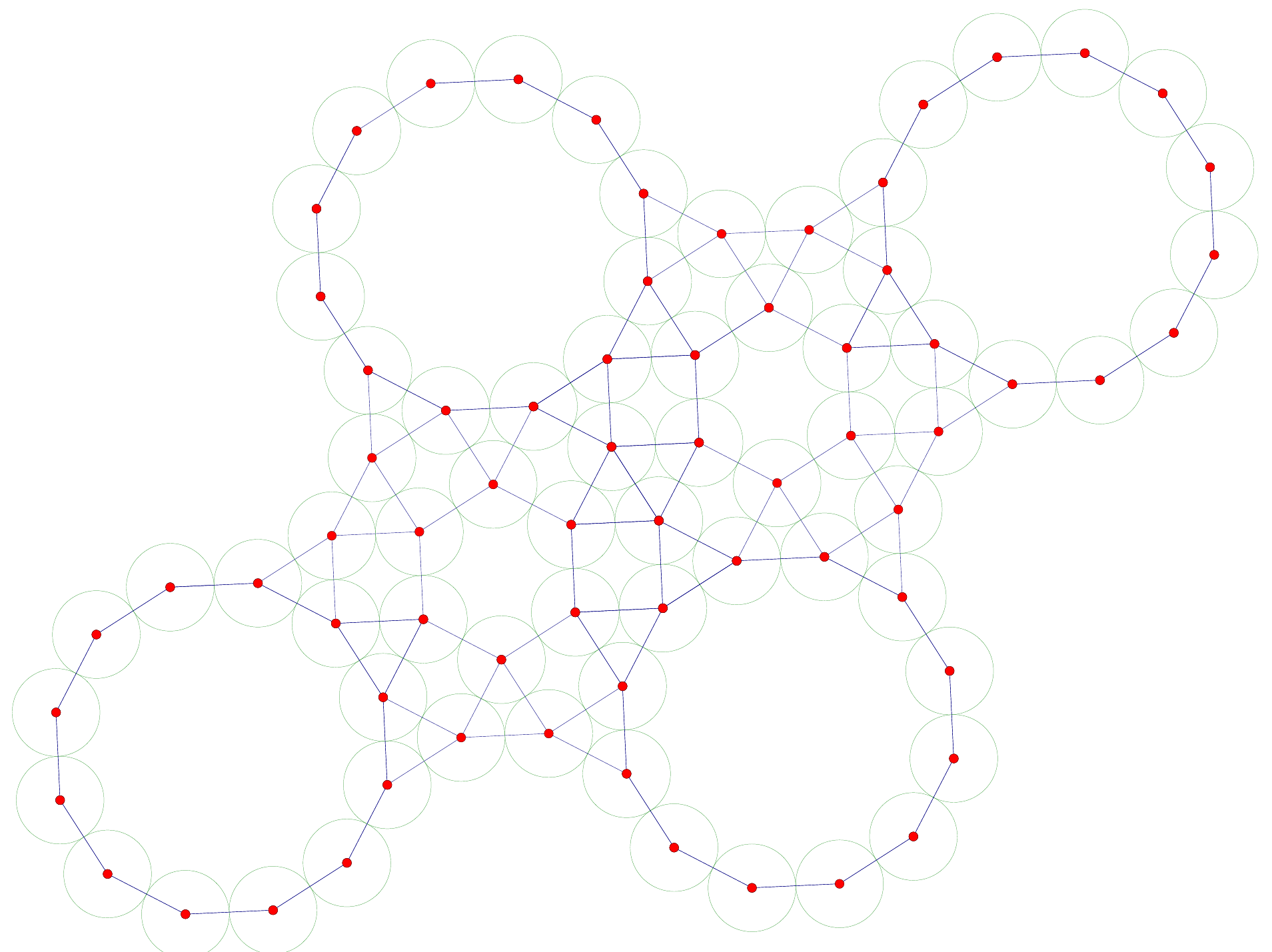}%
        \end{center}
    \caption{This is a packing in a triangular torus with a $12$-gon face that, for a generic configuration of the packing graph is rigid, but it is not rigid as a packing graph since it is a subset of the packing in Figure \ref{fig:triangle-square-final-1}.}
    \label{fig:Ruggero}
    \end{figure}

One of the starting points for this work was from an example of a packing provided by Ruggero Gabbrielli who suggested that the packing in Figure \ref{fig:Ruggero} was at least collectively jammed.
At first it was difficult to determine if, for this packing, it was true, but it was shown by the first author at a talk at the Fields Institute \cite{Connelly-talk} that it was not collectively jammed.  At the next talk by Walter Whiteley \cite{Whiteley-talk} it was shown that by looking at parallel classes of edges of a graph one can create a flex, which is a basic idea of Theorem \ref{thm:main}. It was also shown in \cite{Henley-rotater}, Section 4.9 and Figures 21-22 how the rotating motion described works, and it was shown that there is such a flex when there are two classes of edges as there is here.  

Triangle-square tilings have been of some interest in the physics literature as in \cite{Kawamura, Kawamura-entorpy, Henley, Henley-dodecagon}.  Indeed, the idea of associating to a triangle-square tiling, a graph, where every vertex is of degree six, and identifying it with a triangle lattice, is in \cite{Kawamura}.   In \cite{Kawamura} and \cite{Henley} they estimated the ratio of squares to triangles in a ``random" square-and-triangle tiling.  Square lattice tilings were quite popular.  Suppose that that one side of a period square of square lattice has a polygonal path involving the vectors $1, g^2, g, g^{-1}$, where $g = \sqrt{3}/2  + i/2$ as in Section \ref{sect:TS-tilings}.  So the sum of the edge vectors that go from one corner to the next horizontal corner of the square will be 
\begin{equation*}
s = a + b(g + g^{-1}) = a +b \sqrt{3},
\end{equation*}
where $a$ and $b$ are integers, since $1, g^2, g, g^{-1}$ is another basis for the lattice determined by the edge vectors.  But the area of the square torus is 
\begin{equation*}
s^2 = a^2 + 3b^2 + 2ab \sqrt{3} = f_4 +\frac{\sqrt{3}}{4} f_3 , 
\end{equation*}
where $f_3$ is the number of triangles and $f_4$ is the number of squares.  So 
\begin{equation*}
f_3 = 8 ab\,\,\,\text{and} \,\,f_4 = a^2 + 3 b^2.
\end{equation*}
Thus 
\begin{equation*}
\frac{f_4}{f_3} = \frac{a^2 + 3 b^2}{8 ab} > \frac{1}{8}(\sqrt{3} + \frac{3}{\sqrt{3}})=\frac{\sqrt{3}}{4},
\end{equation*}
since the minimum of $x + 3/x$ for $x > 0$ is when $x = \sqrt{3}$ as any calculus student knows.  

A similar analysis for the triangular torus shows that  $\sqrt{3}/4$ is a maximum for the ratio $f_4/f_3$.

 The ratio $\sqrt{3}/4$ for $f_4/f_3$ was claimed as the maximum entropy for triangle-square tilings with twelve-fold symmetry in \cite{Henley}.
 
\section{Future work}\label{sect:future-work}

It would be helpful to understand the case when the graph of the packing consists of triangles and rhombi more generally than in Theorem \ref{thm:rhombi-triangle}.   There is the question of what kind of packing graphs can occur.   

More ambitiously it would be very useful to understand what kind of collectively jammed packing graphs are that have faces of degree larger than four.  One idea is to complete the graph of the packing to a triangulation, where the degree of each vertex is six.  Then there are ways to compare the density of such packings to the density $\delta_{\Delta}$ and $\frac{n}{n+1}\delta_{\Delta}$ using some simple techniques from the theory of the global rigidity of tensegrities.

\bibliographystyle{plain}
\bibliography{NSF-10,framework}


\end{document}